\documentclass[a4paper,10pt]{article}
\usepackage{stmaryrd}
\usepackage{amsfonts}

\usepackage{bbm}
\usepackage{amscd}
\usepackage{mathrsfs}
\usepackage{latexsym,amssymb,amsmath,amscd,amscd,amsthm,amsxtra}
\usepackage[dvips]{graphicx}
\usepackage[utf8]{inputenc}
\usepackage[T1]{fontenc}
\usepackage{lmodern}
\usepackage{amssymb}
\usepackage[all]{xy}
\usepackage{nicefrac,mathtools,enumitem}
\usepackage{microtype}
\usepackage{xcolor}
\newtheorem{thm}{Theorem}[section]
\newtheorem{lem}[thm]{Lemma}
\newtheorem{cor}[thm]{Corollary}
\newtheorem{pro}[thm]{Proposition}

\newtheorem{defi}[thm]{Definition}

\newcommand {\emptycomment}[1]{}

\newcommand{\lon }{\,\rightarrow\,}
\newcommand{\be }{\begin{equation}}
\newcommand{\ee }{\end{equation}}

\newcommand{\pf}{\noindent{\bf Proof.}\ }

\newcommand{\g}{\frkg}
\newcommand{\h}{\frkh}

\newcommand{\huaB}{\mathcal{B}}%{{\mathcal{E}}}%{\mathcal{B}}

\newcommand{\huaH}{\mathcal{H}}

\newcommand{\huaZ}{\mathcal{Z}}

\newcommand{\frkg}{\mathfrak g}
\newcommand{\frkh}{\mathfrak h}

\def\qed{\hfill ~\vrule height6pt width6pt depth0pt}

\newcommand{\br}[1]{   [ \cdot,    \cdot  ]_\frkg   }

\newcommand{\Id}{\rm{Id}}

\newcommand{\dM}{\mathrm{d}}

\newcommand{\Der}{\mathsf{Der}}
\newcommand{\DER}{\mathsf{DER}}
\newcommand{\Inn}{\mathsf{Inn}}
\newcommand{\Out}{\mathsf{Out}}
\newcommand{\Ad}{\mathsf{Ad}}

\newcommand{\gl}{\mathfrak {gl}}

\newcommand{\cen}{\mathsf{Cen}}

\newcommand{\ad}{\mathsf{ad}}

\newcommand{\K}{\mathbb{K}}

\begin{document}
\title{
{ Cohomology characterizations of non-abelian extensions of  Hom-Lie algebras}
\thanks
 {
Research supported by NSFC (11471139) and NSF of Jilin Province (20170101050JC).
 }
}
\author{Lina Song and Rong Tang  \\
Department of Mathematics, Jilin University,\\
 Changchun 130012, Jilin, China
\\\vspace{3mm}
Email: songln@jlu.edu.cn, tangrong16@mails.jlu.edu.cn}

\date{}


%\footnotetext{The authors declare that there is no conflict of interest regarding the publication of this paper.}

\maketitle
\begin{abstract}
In this paper, first we show that under the assumption of the center of $\h$ being zero, diagonal non-abelian extensions of a regular Hom-Lie algebra $\g$ by a regular Hom-Lie algebra $\h$ are in one-to-one correspondence with  Hom-Lie algebra morphisms from $\g$ to  $\Out(\h)$. Then for a general Hom-Lie algebra morphism from $\g$ to  $\Out(\h)$, we construct a cohomology class as the obstruction of existence of a non-abelian extension that induce the given Hom-Lie algebra morphism.
\end{abstract}

\section{Introduction}

The notion of a Hom-Lie algebra was introduced by Hartwig, Larsson, and Silvestrov in \cite{HLS} as part of a study of deformations of the Witt and the Virasoro algebras. In a Hom-Lie algebra, the Jacobi identity is twisted by a linear map, called the Hom-Jacobi identity. Some $q$-deformations of the Witt and the Virasoro algebras have the structure of a Hom-Lie algebra \cite{HLS,hu}. Because of close relation to discrete and deformed vector fields and differential calculus \cite{HLS,LD1,LD2}, more people pay special attention to this algebraic structure. In particular, representations and deformations of Hom-Lie algebras were studied in \cite{AEM,MS1,sheng}; Extensions of Hom-Lie algebras were studied in \cite{Casas,CasasNonT,LD1}.

 The notion of a Hom-Lie 2-algebra, which is the  categorification of a Hom-Lie algebra, was given in \cite{shengchen}.
In \cite{songtang}, we gave the notion of a derivation of a regular Hom-Lie algebra $(\g,[\cdot,\cdot]_\g,\phi_\g)$. The set of derivations $\Der(\g)$ is a Hom-Lie subalgebra of the regular Hom-Lie algebra $(\gl(\g),[\cdot,\cdot]_{\phi_\g},\Ad_{\phi_\g})$ which was given in \cite{shengxiong}.  We constructed the derivation Hom-Lie 2-algebra $\DER(\g)$, by which we characterize non-abelian extensions of regular Hom-Lie algebras as Hom-Lie 2-algebra morphisms. More precisely, we characterize a diagonal non-abelian extension of a regular Hom-Lie algebra $\g$ by a regular Hom-Lie algebra $\h$ using a Hom-Lie 2-algebra morphism from $\g$ to the derivation Hom-Lie 2-algebra $\DER(\h)$. Associated to a non-abelian extension of a  regular Hom-Lie algebra $\g$ by a regular Hom-Lie algebra $\h$, there is a Hom-Lie algebra morphism from $\g$ to $\Out(\h)$ naturally. However, given an arbitrary Hom-Lie algebra morphism from $\g$ to $\Out(\h)$, whether there is a non-abelian extension of $\g$ by $\h$, that induces the given Hom-Lie algebra morphism, and what is the obstruction is not known yet.

The aim of this paper is to solve the above problem. It turns out that the result is not totally parallel to the case of Lie algebras \cite{AMR,AMR2,Fr,H,IKL}. We need to add some conditions on the short exact sequence related to derivations of Hom-Lie algebras. Under these conditions, first we show that under the assumption of the center of $\h$ being zero, there is a one-to-one correspondence between non-abelian extensions of $\g$ by $\h$ and Hom-Lie algebra morphisms  from $\g$ to $\Out(\h)$. Then for the general case, we show that the obstruction of the existence of a non-abelian extension is given by an element in the third cohomology group.

The paper is organized as follows. In Section 2, we recall some basic notions of Hom-Lie algebras, representations of Hom-Lie algebras, their cohomologies and derivations of  Hom-Lie algebras. In Section 3, we study non-abelian extensions of $\g$ by $\h$ in the case that the center of $\h$ is zero. In Section 4, we give a cohomology characterization of the existence of general  non-abelian extensions of $\g$ by $\h$.

%\vspace{2mm}
 %\noindent {\bf Acknowledgement:} We give our warmest thanks to the referee for very helpful suggestions that improve the paper.

\section{Preliminaries}

 In this paper, we work over an algebraically closed field $\K$ of characteristic 0 and all the vector spaces are over $\K$.
\subsection{Representations, cohomologies and derivations of Hom-Lie algebras}

\begin{defi}\begin{itemize}\item[\rm(1)]
  A $($multiplicative$)$ Hom-Lie algebra is a triple $(\mathfrak{g},[\cdot,\cdot]_{\mathfrak{g}},\phi_{\mathfrak{g}})$ consisting of a
  vector space $\mathfrak{g}$, a skew-symmetric bilinear map (bracket) $[\cdot,\cdot]_{\mathfrak{g}}:\wedge^2\mathfrak{g}\longrightarrow
  \mathfrak{g}$ and a linear map $\phi_{\mathfrak{g}}:\mathfrak{g}\lon \mathfrak{g}$ preserving the bracket, such that  the following Hom-Jacobi
  identity with respect to $\phi_{\mathfrak{g}}$ is satisfied:
  \begin{equation}
    [\phi_{\mathfrak{g}}(x),[y,z]_\g]_\g+[\phi_{\mathfrak{g}}(y),[z,x]_\g]_\g+[\phi_{\mathfrak{g}}(z),[x,y]_\g]_\g=0.
  \end{equation}

\item[\rm(2)]A Hom-Lie algebra is called a regular Hom-Lie algebra if $\phi_{\mathfrak{g}}$ is
an algebra automorphism.
  \end{itemize}
\end{defi}

In the sequel, we always assume that $\phi_{\mathfrak{g}}$ is
an algebra automorphism. That is, in this paper all the Hom-Lie algebras are assumed to be regular Hom-Lie algebras even though some results are also hold for general Hom-Lie algebras.

\begin{defi}
  A morphism of Hom-Lie algebras $f:(\mathfrak{g},[\cdot,\cdot]_{\mathfrak{g}},\phi_{\mathfrak{g}})\lon (\mathfrak{h},[\cdot,\cdot]_{\mathfrak{h}},\phi_{\mathfrak{h}})$ is a linear map $f:\mathfrak{g}\lon \mathfrak{h}$ such that
  \begin{eqnarray}
f[x,y]_{\mathfrak{g}}&=&[f(x),f(y)]_{\mathfrak{h}},\hspace{3mm}\forall x,y\in \mathfrak{g},\\
    f\circ \phi_{\mathfrak{g}}&=&\phi_{\mathfrak{h}}\circ f.
  \end{eqnarray}
\end{defi}

\begin{defi}
A representation of a Hom-Lie algebra $(\mathfrak{g},[\cdot,\cdot]_{\frkg},\phi_{\mathfrak{g}})$ on a vector space $V$ with respect to $\beta \in\mathfrak{gl}(V)$ is a linear map $\rho:\mathfrak{g}\lon \mathfrak{gl}(V)$ such that for all $x,y\in\mathfrak{g}$, the following equalities are satisfied:
\begin{eqnarray}
\rho(\phi_{\mathfrak{g}}(x))\circ \beta&=&\beta\circ\rho(x),\\
\rho([x,y]_{\mathfrak{g}})\circ\beta&=&\rho(\phi_{\mathfrak{g}}(x))\circ\rho(y)-\rho(\phi_{\mathfrak{g}}(y))\circ\rho(x).
\end{eqnarray}
\end{defi}

We denote a representation by $(\rho,V,\beta)$. For all $x\in\mathfrak{g}$, we define $\ad_{x}:\mathfrak{g}\lon \mathfrak{g}$ by
\begin{eqnarray}
\ad_{x}(y)=[x,y]_{\mathfrak{g}},\quad\forall y \in \mathfrak{g}.
\end{eqnarray}
Then $\ad:\g\longrightarrow\gl(\frak g)$ is a representation of the Hom-Lie algebra $(\mathfrak{g},[\cdot,\cdot]_{\mathfrak{g}},\phi_{\mathfrak{g}})$ on $\g$ with respect to $\phi_\g$, which is called the {\bf adjoint representation}.

Let  $(\rho,V,\beta)$ be a representation. We define the set of $k$-Hom-cochains by
 \begin{eqnarray*}
 C_{\phi_{\g},\beta}^{k}(\g;V)\stackrel{\triangle}{=}\{f\in C^{k}(\g;V)\mid \beta\circ \phi_{\g}=\phi_{\g}\circ\beta^{\otimes k}\}.
 \end{eqnarray*}
For $k\geq1$, we define the coboundary operator $\dM_{\rho}:C_{\phi_{\g},\beta}^{k}(\g;V)\lon C_{\phi_{\g},\beta}^{k+1}(\g;V)$ by
\begin{eqnarray*}
(\dM_{\rho}f)(x_1,\cdots,x_{k+1})&=&\sum_{i=1}^{k+1}(-1)^{i+1}\rho(\phi_{\g}^{k-1}(x_i))(f(x_1,\cdots,\widehat{x_{i}},\cdots,x_{k+1}))\\
                         &&+\sum_{i<j}(-1)^{i+j}f([x_i,x_j]_{\g},\phi_{\g}(x_1),\cdots,\widehat{\phi_{\g}(x_i)},\cdots,\widehat{\phi_{\g}(x_j)},\cdots,\phi_{\g}(x_{k+1})).
\end{eqnarray*}
The fact that $\dM_{\rho}\circ\dM_{\rho}=0$ is proved in \cite{sheng}.  Denote by $\mathcal{Z}^k(\g;\rho)$ and $\mathcal{B}^k(\g;\rho)$ the sets of $k$-cocycles and  $k$-coboundaries respectively. We define the $k$-th cohomology group
$\mathcal{H}^k(\g;\rho)$ to be $\mathcal{Z}^k(\g;\rho)/\mathcal{B}^k(\g;\rho)$.

Let $V$ be a vector space, and $\beta\in GL(V)$. Define a skew-symmetric bilinear bracket operation $[\cdot,\cdot]_{\beta}:\wedge^2\mathfrak{gl}(V)\longrightarrow\mathfrak{gl}(V)$ by
\begin{eqnarray}\label{eq:bracket}
[A,B]_{\beta}=\beta \circ A \circ\beta^{-1}\circ B \circ\beta^{-1}-\beta\circ B \circ\beta^{-1}\circ A\circ \beta^{-1}, \hspace{3mm}\forall A,B\in \mathfrak{gl}(V).
\end{eqnarray}
Denote by $\Ad_{\beta}:\mathfrak{gl}(V)\lon \mathfrak{gl}(V)$ the adjoint action on $\mathfrak{gl}(V)$, i.e.
\begin{equation}\label{eq:Ad}
\Ad_{\beta}(A)=\beta\circ A\circ \beta^{-1}.
\end{equation}
Then $(\mathfrak{gl}(V),[\cdot,\cdot]_{\beta},\Ad_{\beta})$ is a regular Hom-Lie algebra. See \cite{shengxiong}
for more details.

\begin{thm}{\rm (\cite[Theorem 4.2]{shengxiong})}\label{thm:Hom-Lie rep}
Let $(\g,[\cdot,\cdot]_\g,\phi_\g)$ be a Hom-Lie algebra, $V$ a vector space, $\beta \in GL(V)$. Then, $\rho:\g\lon\gl(V)$ is a representation of $(\g,[\cdot,\cdot]_\g,\phi_\g)$ on $V$ with respect to $\beta$ if and only if $\rho:(\g,[\cdot,\cdot]_\g,\phi_\g)\lon
(\mathfrak{gl}(V),[\cdot,\cdot]_{\beta},\Ad_{\beta})$ is a morphism of Hom-Lie algebras.
\end{thm}

\begin{defi}{\rm (\cite[Definition 3.1]{songtang})}\label{def:derivation}
A linear map $D:\mathfrak{g}\lon \mathfrak{g}$ is called a derivation of a Hom-Lie algebra $(\mathfrak{g},[\cdot,\cdot]_{\mathfrak{g}},\phi_{\mathfrak{g}})$ if
\begin{eqnarray}\label{derivation}
D[x,y]_{\mathfrak{g}}=[\phi_{\mathfrak{g}}(x),(\Ad_{\phi_{\mathfrak{g}}^{-1}}D)(y)]_{\mathfrak{g}}+[(\Ad_{\phi_{\mathfrak{g}}^{-1}}D)(x),
\phi_{\mathfrak{g}}(y)]_{\mathfrak{g}},\hspace{3mm}\forall x,y\in\mathfrak{g}.
\end{eqnarray}
\end{defi}
Denote by $\Der(\mathfrak{g})$ the set of derivations of the Hom-Lie algebra $(\mathfrak{g},[\cdot,\cdot]_{\mathfrak{g}},\phi_{\mathfrak{g}})$.
Then we obtain that $(\Der(\g),[\cdot,\cdot]_{\phi_{\frkg}},\Ad_{\phi_{\mathfrak{g}}})$ is a Hom-Lie algebra, which is a subalgebra of the Hom-Lie algebra $(\mathfrak{gl}(\frkg),[\cdot,\cdot]_{\phi_{\frkg}},\Ad_{\phi_{\frkg}})$.

 For all $x\in\g$, $\ad_{x}$ is a derivation of the Hom-Lie algebra $(\mathfrak{g},[\cdot,\cdot]_{\mathfrak{g}},\phi_{\mathfrak{g}})$, which we call an {\bf inner derivation}. See \cite{songtang} for more details.
Denote by $\Inn(\g)$ the set of inner derivations of the Hom-Lie algebra $(\mathfrak{g},[\cdot,\cdot]_{\mathfrak{g}},\phi_{\mathfrak{g}})$, i.e.
\begin{eqnarray}
\Inn(\mathfrak{g})=\{\ad_{x}\mid x\in \mathfrak{g}\}.
\end{eqnarray}

\begin{lem}{\rm (\cite[Lemma 3.6]{songtang})}\label{lem:ideal}
Let $(\mathfrak{g},[\cdot,\cdot]_{\mathfrak{g}},\phi_{\mathfrak{g}})$ be a Hom-Lie algebra.  For all $x\in\g$ and $D\in\Der(\g)$, we have
  $$\Ad_{\phi_\g}\ad_x=\ad_{\phi_{\g}(x)},\quad [D,\ad_x]_{\phi_\g}=\ad_{D(x)}.$$
  Therefore, $\Inn(\g)$ is an ideal of the Hom-Lie algebra $(\Der(\g),[\cdot,\cdot]_{\phi_{\frkg}},\Ad_{\phi_{\mathfrak{g}}})$.
\end{lem}

 Denote by $\Out(\g)$ the set of out derivations of the Hom-Lie algebra $(\mathfrak{g},[\cdot,\cdot]_{\mathfrak{g}},\phi_{\mathfrak{g}})$, i.e.
\begin{eqnarray}
\Out(\g)=\Der(\g)/\Inn(\g).
\end{eqnarray}
We use $\pi$ to denote the quotient map from $\Der(\g)$ to $\Out(\g)$.
\subsection{Non-abelian extensions of Hom-Lie algebras}

\begin{defi}
A non-abelian extension of a Hom-Lie algebra $(\mathfrak{g},[\cdot,\cdot]_{\mathfrak{g}},\phi_{\mathfrak{g}})$ by a Hom-Lie algebra $(\mathfrak{h},[\cdot,\cdot]_{\mathfrak{h}},\phi_{\mathfrak{h}})$ is a commutative diagram with rows being short exact sequence of Hom-Lie algebra morphisms:
\[\begin{CD}
0@>>>\mathfrak{h}@>\iota>>\hat{\mathfrak{g}}@>p>>{\mathfrak{g}}             @>>>0\\
@.    @V\phi_{\mathfrak{h}}VV   @V\phi_{\hat{\mathfrak{g}}}VV  @V\phi_{\mathfrak{g}}VV    @.\\
0@>>>\mathfrak{h}@>\iota>>\hat{\mathfrak{g}}@>p>>{\mathfrak{g}}             @>>>0
,\end{CD}\]
where $(\hat{\mathfrak{g}},[\cdot,\cdot]_{\hat{\g}},\phi_{\hat{\mathfrak{g}}})$ is a Hom-Lie algebra.
\end{defi}

We can regard $\mathfrak{h}$ as a subspace of $\hat{\mathfrak{g}}$ and $\phi_{\hat{\mathfrak{g}}}|_{\mathfrak{h}}=\phi_{\mathfrak{h}}$. Thus, $\mathfrak{h}$ is an invariant subspace of $\phi_{\hat{\mathfrak{g}}}$.  We say that an extension is {\bf diagonal} if
$\hat{\mathfrak{g}}$ has an invariant subspace $X$ of $\phi_{\hat{\mathfrak{g}}}$ such that $\mathfrak{h}\oplus X=\hat{\mathfrak{g}}$. In general, $\hat{\mathfrak{g}}$ does not always have an invariant subspace $X$ of $\phi_{\hat{\mathfrak{g}}}$ such that $\mathfrak{h}\oplus X=\hat{\mathfrak{g}}$. For example, the matrix representation of $\phi_{\hat{\mathfrak{g}}}$ is a Jordan block. We only study diagonal non-abelian extensions in the sequel.

\begin{defi}\label{defi:iso}
Two extensions of $\mathfrak{g}$ by $\mathfrak{h}$, $(\hat{\g}_1,[\cdot,\cdot]_{\hat{\g}_1},\phi_{\hat{\g}_1})$ and $(\hat{\g}_2,[\cdot,\cdot]_{\hat{\g}_2},\phi_{\hat{\g}_2})$, are said to be isomorphic if there exists a Hom-Lie algebra morphism $\theta:\hat{\mathfrak{g}}_{2}\lon \hat{\mathfrak{g}}_{1}$ such that we have the following commutative diagram:
\label{seqiso}\[\begin{CD}
0@>>>\mathfrak{h}@>\iota_{2}>>\hat{\mathfrak{g}}_{2}@>p_{2}>>{\mathfrak{g}}             @>>>0\\
@.    @|                       @V\theta VV                     @|                       @.\\
0@>>>\mathfrak{h}@>\iota_{1}>>\hat{\mathfrak{g}}_{1}@>p_{1}>>{\mathfrak{g}}             @>>>0
.\end{CD}\]
\end{defi}

\begin{lem}{\rm (\cite[Lemma 4.4]{songtang})}
A Hom-Lie algebra $(\hat{\g},[\cdot,\cdot]_{\hat{\g}},\phi_{\hat{\g}})$ is a diagonal non-abelian extension of a Hom-Lie algebra $(\g,[\cdot,\cdot]_{\g},\phi_\g)$ by a Hom-Lie algebra $(\h,[\cdot,\cdot]_{\h},\phi_\h)$ if and only if there is a section $s:\mathfrak{g}\lon \hat{\mathfrak{g}}$  such that
$p\circ s=\Id$ and $\phi_{\hat{\mathfrak{g}}}\circ s=s\circ\phi_{\mathfrak{g}}$. This section is called {\bf diagonal}.
\end{lem}

Let $(\hat{\g},[\cdot,\cdot]_{\hat{\g}},\phi_{\hat{\g}})$ be a diagonal extension of a Hom-Lie algebra $(\g,[\cdot,\cdot]_{\g},\phi_\g)$ by a Hom-Lie algebra $(\h,[\cdot,\cdot]_{\h},\phi_\h)$ and $s:\frkg\lon \hat{\frkg}$ a diagonal section.
Define linear maps $\omega:\mathfrak{g}\wedge\mathfrak{g}\lon \mathfrak{h}$ and $\rho:\frkg\lon \mathfrak{gl}(\mathfrak{h})$   respectively by
\begin{eqnarray}
\label{do}\omega(x,y)&=&[s(x),s(y)]_{\hat{\mathfrak{g}}}-s[x,y]_{\mathfrak{g}}, \,\,\,\,\forall x,y \in \mathfrak{g},\\
\label{dr}\rho_{x}(u)&=&[s(x),u]_{\hat{\mathfrak{g}}}, \,\,\,\,\forall x\in \mathfrak{g}, u\in \mathfrak{h}.
\end{eqnarray}
Obviously, $\hat{\g}$ is isomorphic to $\g\oplus\h$ as vector spaces. Transfer the Hom-Lie algebra structure on $\hat{\mathfrak{g}}$ to that on $\mathfrak{g}\oplus \mathfrak{h}$, we obtain a Hom-Lie algebra $(\mathfrak{g}\oplus \mathfrak{h},[\cdot,\cdot]_{(\rho,\omega)},\phi)$, where $[\cdot,\cdot]_{(\rho,\omega)}$ and $\phi$ are given by
\begin{eqnarray}
\label{dbr}[x+u,y+v]_{(\rho,\omega)}&=&[x,y]_{\mathfrak{g}}+\omega(x,y)+\rho_{x}(v)-\rho_{y}(u)+[u,v]_{\mathfrak{h}},\\
 \label{dmo}\phi(x+u)&=&\phi_{\mathfrak{g}}(x)+\phi_{\mathfrak{h}}(u).
\end{eqnarray}
The following proposition gives the conditions on $\rho$ and $\omega$ such that $(\g\oplus\h, [\cdot,\cdot]_{(\rho,\omega)},\phi)$ is a Hom-Lie algebra.

\begin{pro}{\rm (\cite[Proposition 4.5]{songtang})}\label{extension}
With the above notations, $(\g\oplus\h, [\cdot,\cdot]_{(\rho,\omega)},\phi)$ is a Hom-Lie algebra if and only if $\rho$ and $\omega $ satisfy the following equalities:
\begin{eqnarray}
\label{p1}\phi_{\mathfrak{h}}\circ \rho_{x}&=&\rho_{\phi_{\mathfrak{g}}(x)}\circ \phi_{\mathfrak{h}},\\
\label{p2}\rho_{x}([u,v]_{\mathfrak{h}})&=&[\phi_{\mathfrak{h}}(u),(\Ad_{\phi_{\mathfrak{h}}^{-1}}\rho_{x})(v)]_{\mathfrak{h}}+
[(\Ad_{\phi_{\mathfrak{h}}^{-1}}\rho_{x})(u),\phi_{\mathfrak{h}}(v)]_{\mathfrak{h}},\\[-6pt]
\label{p3}\phi_{\mathfrak{h}}\circ \omega&=&\omega\circ\phi_{\mathfrak{g}}^{\otimes 2},\\[-6pt]
\label{p4}[\rho_{x},\rho_{y}]_{\phi_{\mathfrak{h}}}-\rho_{[x,y]_{\mathfrak{g}}}&=&\ad_{\omega(x,y)},\\
\label{p5}\rho_{\phi_{\mathfrak{g}}(x)}(\omega(y,z))+c.p.&=&\omega([x,y]_{\mathfrak{g}},\phi_{\mathfrak{g}}(z))+c.p..
\end{eqnarray}
\end{pro}

 For any diagonal non-abelian extension, by choosing a diagonal section, it is isomorphic to
 $(\mathfrak{g}\oplus \mathfrak{h},[\cdot,\cdot]_{(\rho,\omega)},\phi)$. Therefore, we only consider diagonal non-abelian extensions of the form $(\mathfrak{g}\oplus \mathfrak{h},[\cdot,\cdot]_{(\rho,\omega)},\phi)$ in the sequel.

\begin{pro}{\rm (\cite[Theorem 4.9]{songtang})}\label{iso}
Let $(\mathfrak{g}\oplus \mathfrak{h},[\cdot,\cdot]_{(\rho,\omega)},\phi)$ and $(\mathfrak{g}\oplus \mathfrak{h},[\cdot,\cdot]_{(\rho',\omega')},\phi)$ be two diagonal non-abelian extensions of $\mathfrak{g}$ by $\mathfrak{h}$. The two extensions are equivalent if and only if there is a linear map $\xi:\g\lon\h$ such that
\begin{eqnarray}
\label{isom1}\phi_{\mathfrak{h}}(\xi(x))&=&\xi(\phi_{\mathfrak{g}}(x)),\\
\label{isom2}\rho'_{x}-\rho_{x}&=&\ad_{\xi(x)},\\
\label{isom3}\omega'(x,y)-\omega(x,y)&=&\rho_{x}(\xi(y))-\rho_{y}(\xi(x))+[\xi(x),\xi(y)]_{\mathfrak{h}}
-\xi([x,y]_{\mathfrak{g}}).
\end{eqnarray}
\end{pro}

\section{Classification of diagonal non-abelian extensions of Hom-Lie algebras: special case}
In this section, we classify diagonal non-abelian extensions of Hom-Lie algebras for the case that $\cen(\h)=0.$

\begin{thm}\label{thm:special}
Let $(\g,[\cdot,\cdot]_\g,\phi_\g)$ and $(\h,[\cdot,\cdot]_\h,\phi_\h)$ be Hom-Lie algebras such that $\cen(\h)=0$. If the following short exact sequence of Hom-Lie algebra morphisms
\begin{equation}\label{seq:non-abelianext}
0\lon\Inn(\h)\stackrel{}{\lon}\Der(\h)\stackrel{\pi}{\lon} \Out(\h)\lon 0,
\end{equation}
is a diagonal non-abelian extension of $\Out(\h)$ by $\Inn(\h)$,
then isomorphism classes of diagonal non-abelian extensions of $\g$ by $\h$ correspond bijectively to Hom-Lie algebra homomorphisms
\begin{eqnarray*}
\bar{\rho}:\g\lon \Out(\h).
\end{eqnarray*}
\end{thm}

\pf Let $(\mathfrak{g}\oplus \mathfrak{h},[\cdot,\cdot]_{(\rho,\omega)},\phi)$ be a diagonal non-abelian extensions of $\g$ by $\h$ given by
\eqref{dbr}-\eqref{dmo}. By \eqref{p2}, we have $\rho_x\in\Der(\h)$. Let $\pi:\Der(\frkh)\lon \Out(\frkh)$ be the quotient map. We denote  the induced Hom-Lie algebra structure on $\Out(\h)$ by $[\cdot,\cdot]'_{\phi_\h}$ and $\Ad_{\phi_\h}'$. Hence we can define
$$\bar{\rho}=\pi\circ\rho.$$
By \eqref{p1}, for all $x\in\g$ we have
$$\bar{\rho}_{\phi_\g(x)}=\pi(\rho_{\phi_{\g}(x)})=\pi(\Ad_{\phi_\h}(\rho_x))=\Ad_{\phi_\h}'(\bar{\rho}_x).$$
By \eqref{p4}, we have
$$\bar{\rho}_{[x,y]_\g}=\pi([\rho_x,\rho_y]_{\phi_\h}-\ad_{\omega(x,y)})=\pi([\rho_x,\rho_y]_{\phi_\h})=[\bar{\rho}_x,\bar{\rho}_y]_{\phi_\mathfrak{h}}'.$$
Thus, $\bar{\rho}$ is a Hom-Lie algebra homomorphism from $\g$ to $\Out(\h)$.

Let $(\mathfrak{g}\oplus \mathfrak{h},[\cdot,\cdot]_{(\rho',\omega')},\phi)$ and $(\mathfrak{g}\oplus \mathfrak{h},[\cdot,\cdot]_{(\rho,\omega)},\phi)$ be isomorphic diagonal non-abelian extensions of $\mathfrak{g}$ by $\mathfrak{h}$. By Proposition \ref{iso}, we have
$$\bar{\rho}'_x=\pi(\rho_x+\ad_{\xi(x)})=\bar{\rho}_x.$$
Thus, we obtain that isomorphic diagonal non-abelian extensions of $\g$ by $\h$ correspond to the same Hom-Lie algebra homomorphism from $\g$ to $\Out(\h)$.

Conversely, let $\bar{\rho}$ be a Hom-Lie algebra homomorphism from $\g$ to $\Out(\h)$. Since the short exact sequence of Hom-Lie algebras \eqref{seq:non-abelianext} is a diagonal non-abelian extensions of $\Out(\h)$ by $\Inn(\h)$, we can choose a diagonal section $s$ of $\pi:\Der(\g)\lon\Out(\h)$. Moreover, we define $\rho:\g\lon\gl(\h)$ by
\begin{eqnarray}\label{rho}
\rho_x=s(\bar{\rho}_x).
\end{eqnarray}
We have $\rho_x\in\Der(\h)$. Thus we get \eqref{p2}. Since $s$ is a diagonal section, we have
\begin{eqnarray}\nonumber
\rho_{\phi_{\g}(x)}\circ \phi_\mathfrak{h}&=&\big(s(\bar{\rho}_{\phi_{\g}(x)})\big)\circ\phi_\mathfrak{h}=\big(s(\Ad_{\phi_\mathfrak{h}}'(\bar{\rho}_x))\big)\circ\phi_\mathfrak{h}
                                                              =\big((s\circ \Ad_{\phi_\mathfrak{h}}')\circ\bar{\rho}_x)\big)\circ\phi_\mathfrak{h}\\\nonumber
                                                              &=&\big((\Ad_{\phi_\mathfrak{h}}\circ s)\circ\bar{\rho}_x)\big)\circ\phi_\mathfrak{h}
                                                              =\big(\Ad_{\phi_\mathfrak{h}}(s(\bar{\rho}_x))\big)\circ\phi_\mathfrak{h}
                                                              =\phi_\mathfrak{h}\circ s(\bar{\rho}_x)\circ \phi_\mathfrak{h}^{-1} \circ\phi_\mathfrak{h}\\
                                                             &=&\label{ext1}\phi_\mathfrak{h}\circ \rho_x.
\end{eqnarray}
Thus, we obtain \eqref{p1}. Since $\pi$ and $\bar{\rho}$ are Hom-Lie algebra homomorphisms, for all $x,y\in\g$, we have
$$\pi([\rho_x,\rho_y]_{\phi_\mathfrak{h}}-\rho_{[x,y]_\g})=[\pi(\rho_x),\pi(\rho_y)]_{\phi_\mathfrak{h}}'-\pi(\rho_{[x,y]_\g})=[\bar{\rho}_x,\bar{\rho}_y]_{\phi_\mathfrak{h}}'-
\bar{\rho}_{[x,y]_\g}=0,$$
which implies that $[\rho_x,\rho_y]_{\phi_\mathfrak{h}}-\rho_{[x,y]_\g}\in\Inn(\h)$. Since we have the following short exact sequence of Hom-Lie algebra morphisms
\begin{equation*}\label{seq:abelianext}
0\lon \cen(\h)\stackrel{}{\lon}\mathfrak{h}\stackrel{\ad}{\lon} \Inn(\mathfrak{h})\lon 0,
\end{equation*}
and $\cen(\h)=0$, there exists a unique linear map $\omega:\g\wedge\g\lon\h$ such that
\begin{eqnarray}
\label{ext3}[\rho_{x},\rho_{y}]_{\phi_{\mathfrak{h}}}-\rho_{[x,y]_{\g}}&=&\ad_{\omega(x,y)}.
\end{eqnarray}
Furthermore, we claim that
\begin{equation}
  \label{ext2}\phi_{\mathfrak{h}}\circ \omega=\omega\circ\phi_{\mathfrak{g}}^{\otimes 2}.
\end{equation}
In fact, for all $u\in\h$, we have
\begin{eqnarray*}
  &&[\phi_{\mathfrak{h}}( \omega(x,y))-\omega(\phi_{\mathfrak{g}}(x),\phi_\g(y)),\phi_\h(u)]_\h\\&=&\phi_\h\ad_{\omega(x,y)}u-\ad_{\omega(\phi_{\mathfrak{g}}(x),\phi_\g(y))}\phi_\h(u)\\
 &=& \phi_\h([\rho_{x},\rho_{y}]_{\phi_{\mathfrak{h}}}(u)-\rho_{[x,y]_{\g}}(u))
 -\big([\rho_{\phi_\g(x)},\rho_{\phi_\g(y)}]_{\phi_{\mathfrak{h}}}(\phi_\h(u))-\rho_{\phi_\g[x,y]_{\g}}(\phi_\h(u))\big)\\
 &=&0,
\end{eqnarray*}
which implies that $\phi_{\mathfrak{h}}( \omega(x,y))-\omega(\phi_{\mathfrak{g}}(x),\phi_\g(y))=0$ since $\cen(\h)=0.$
Thus, we obtain \eqref{p3} and \eqref{p4}. For all $x,y,z\in \g,u\in \mathfrak{h}$, by $\rho_x\in\Der(\h)$ and \eqref{ext1}-\eqref{ext3}, we have
\begin{eqnarray*}
&&\,\,\,\,[\big(\rho_{\phi_{\mathfrak{g}}(x)}(\omega(y,z))-\omega([x,y]_{\mathfrak{g}},\phi_{\mathfrak{g}}(z))\big)+c.p.,\phi_\mathfrak{h}(u)]_\mathfrak{h}\\
%&&=\phi_\mathfrak{h}[\big((\phi_\mathfrak{h}^{-1}\circ\rho_{\phi_{\mathfrak{g}}(x)}\circ\phi_\mathfrak{h})(\omega(\phi_\g^{-1}y,\phi_\g^{-1}z))-
%\omega([\phi_\g^{-1}x,\phi_\g^{-1}y]_{\mathfrak{g}},z)\big)+c.p.,u]_\mathfrak{h}\\
&&=\phi_\mathfrak{h}\big([(\Ad_{\phi_{\mathfrak{h}}^{-1}}\rho_{\phi_{\g}(x)})(\omega(\phi_\g^{-1}y,\phi_\g^{-1}z)),\phi_{\mathfrak{h}}(\phi_{\mathfrak{h}}^{-1}u)]_\mathfrak{h}-
[\omega([\phi_\g^{-1}x,\phi_\g^{-1}y]_{\mathfrak{g}},z),u]_\mathfrak{h}+c.p.\big)\\
&&=\phi_\mathfrak{h}\big(\rho_{\phi_{\g}(x)}[\omega(\phi_\g^{-1}y,\phi_\g^{-1}z),\phi_{\mathfrak{h}}^{-1}u]_\mathfrak{h}-
[\omega(y,z),(\Ad_{\phi_{\mathfrak{h}}^{-1}}\rho_{\phi_{\g}(x)})(\phi_{\mathfrak{h}}^{-1}u)]_\mathfrak{h}\\
&&\,\,\,\,\,\,\,-[\rho_{[\phi_{\g}^{-1}x,\phi_{\g}^{-1}y]_{\g}},\rho_{z}]_{\phi_{\mathfrak{h}}}(u)+\rho_{[[\phi_{\g}^{-1}x,\phi_{\g}^{-1}y]_{\g},z]_{\g}}(u)+c.p.\big)\\
&&=\phi_\mathfrak{h}\big(\rho_{\phi_{\g}(x)}\big(([\rho_{\phi_{\g}^{-1}y},\rho_{\phi_{\g}^{-1}z}]_{\phi_{\mathfrak{h}}}-\rho_{[\phi_{\g}^{-1}y,\phi_{\g}^{-1}z]_{\g}})(\phi_{\mathfrak{h}}^{-1}u)\big)               -([\rho_{y},\rho_{z}]_{\phi_{\mathfrak{h}}}-\rho_{[y,z]_{\g}})((\rho_{x}\circ\phi_{\mathfrak{h}}^{-1})(u))\\
&&\,\,\,\,\,\,\,-[\rho_{[\phi_{\g}^{-1}x,\phi_{\g}^{-1}y]_{\g}},\rho_{z}]_{\phi_{\mathfrak{h}}}(u)+\rho_{[[\phi_{\g}^{-1}x,\phi_{\g}^{-1}y]_{\g},z]_{\g}}(u)+c.p.\big)\\
&&=\phi_{\mathfrak{h}}\big((\phi_{\mathfrak{h}}\circ\rho_{x}\circ\phi_{\mathfrak{h}}^{-1}\circ [\rho_{\phi_{\g}^{-1}y},\rho_{\phi_{\g}^{-1}z}]_{\phi_{\mathfrak{h}}}\circ\phi_{\mathfrak{h}}^{-1}-\phi_{\mathfrak{h}}\circ[\rho_{\phi_{\g}^{-1}y},\rho_{\phi_{\g}^{-1}z}]_{\phi_{\mathfrak{h}}}\circ\phi_{\mathfrak{h}}^{-1}\circ\rho_{x}\circ\phi_{\mathfrak{h}}^{-1})(u)\\
&&\,\,\,\,\,\,\,+(\phi_{\mathfrak{h}}\circ\rho_{[\phi_{\g}^{-1}y,\phi_{\g}^{-1}z]_{\g}}\circ\phi_{\mathfrak{h}}^{-1}\circ\rho_{x}\circ\phi_{\mathfrak{h}}^{-1}-\phi_{\mathfrak{h}}\circ\rho_{x}\circ\phi_{\mathfrak{h}}^{-1}\circ\rho_{[\phi_{\g}^{-1}y,\phi_{\g}^{-1}z]_{\g}}\circ\phi_{\mathfrak{h}}^{-1})(u)\\
&&\,\,\,\,\,\,\,-[\rho_{[\phi_{\g}^{-1}x,\phi_{\g}^{-1}y]_{\g}},\rho_{z}]_{\phi_{\mathfrak{h}}}(u)+\rho_{[[\phi_{\g}^{-1}x,\phi_{\g}^{-1}y]_{\g},z]_{\g}}(u)+c.p.
\big)\\
&&=\phi_{\mathfrak{h}}\big(([\rho_{x},[\rho_{\phi_{\g}^{-1}y},\rho_{\phi_{\g}^{-1}z}]_{\phi_{\mathfrak{h}}}]_{\phi_{\mathfrak{h}}}+[\rho_{[\phi_{\g}^{-1}y,\phi_{\g}^{-1}z]_{\g}},\rho_{x}]_{\phi_{\mathfrak{h}}}\\
&&\,\,\,\,\,\,\,-[\rho_{[\phi_{\g}^{-1}x,\phi_{\g}^{-1}y]_{\g}},\rho_{z}]_{\phi_{\mathfrak{h}}}+\rho_{[[\phi_{\g}^{-1}x,\phi_{\g}^{-1}y]_{\g},z]_{\g}})(u)+c.p.\big)\\
&&=0.
\end{eqnarray*}
Thus, we have
\begin{eqnarray}\label{cen}
\rho_{\phi_{\mathfrak{g}}(x)}(\omega(y,z))-\omega([x,y]_{\mathfrak{g}},\phi_{\mathfrak{g}}(z))+c.p.\in\cen(\h).
\end{eqnarray}
Since $\cen(\h)=0$, we have \eqref{p5}. Therefore, we deduce that \eqref{p1}-\eqref{p5} hold. By Proposition \ref{extension},
$(\g\oplus\h, [\cdot,\cdot]_{(\rho,\omega)},\phi)$ is a diagonal non-abelian extension of $\mathfrak{g}$ by $\mathfrak{h}$.

If we choose another section $s'$ of  $\pi:\Der(\g)\lon\Out(\h)$, we obtain another diagonal non-abelian extension $(\mathfrak{g}\oplus \mathfrak{h},[\cdot,\cdot]_{(\rho',\omega')},\phi)$. Obviously, we have
\begin{eqnarray*}
\pi(\rho_{x}'-\rho_{x})=(\pi\circ s')(\bar{\rho}_x)-(\pi\circ s)(\bar{\rho}_x)=0,
\end{eqnarray*}
which implies that $\rho_{x}'-\rho_{x}\in \Inn(\mathfrak{h})$. Since $\cen(\h) = 0$, there is a unique linear map
$\xi:\g\lon\mathfrak{h}$ such that
\begin{eqnarray}
\label{ext-isom2}\rho'_{x}-\rho_{x}&=&\ad_{\xi(x)}.
\end{eqnarray}
We claim that
\begin{equation}
  \label{ext-isom1}\phi_{\mathfrak{h}}(\xi(x))=\xi(\phi_{\mathfrak{g}}(x)).
\end{equation}
In fact, for all $u\in\h$, we have
\begin{eqnarray*}
  [\phi_{\mathfrak{h}}(\xi(x))-\xi(\phi_{\mathfrak{g}}(x)),\phi_\h(u)]_\h&=&\phi_\h\ad_{\xi(x)}u-\ad_{\xi(\phi_\g(x))}\phi_\h(u)\\
  &=&\phi_\h(\rho'_xu-\rho_xu)-(\rho'_{\phi_\g(x)}\phi_\h(u)-\rho_{\phi_\g(x)}\phi_\h(u))\\
  &=&0,
\end{eqnarray*}
which implies that $\phi_{\mathfrak{h}}(\xi(x))=\xi(\phi_{\mathfrak{g}}(x))$ since $\cen(\h)=0$. Thus, we obtain \eqref{isom1} and \eqref{isom2}.
By Lemma \ref{lem:ideal} and \eqref{ext3}, for all $x, y\in\g$ we have
\begin{eqnarray*}
\ad_{\omega'(x,y)-\omega(x,y)}&=&[\rho'_x,\rho'_y]_{\phi_\mathfrak{h}}-\rho'_{[x,y]_\g}-[\rho_x,\rho_y]_{\phi_\mathfrak{h}}+\rho_{[x,y]_\g}\\
                           &=&[\rho_x+\ad_{\xi(x)},\rho_y+\ad_{\xi(y)}]_{\phi_\mathfrak{h}}-\rho_{[x,y]_\g}-\ad_{\xi([x,y]_\g)}
                           -[\rho_x,\rho_y]_{\phi_\mathfrak{h}}+\rho_{[x,y]_\g}\\
                           &=&[\rho_x,\ad_{\xi(y)}]_{\phi_\mathfrak{h}}+[\ad_{\xi(x)},\rho_y]_{\phi_\mathfrak{h}}+
                   [\ad_{\xi(x)},\ad_{\xi(y)}]_{\phi_\mathfrak{h}}-\ad_{\xi([x,y]_\g)}\\
                           &=&\ad_{\rho_x(\xi(y))}-\ad_{\rho_y(\xi(x))}+\ad_{[\xi(x),\xi(y)]_{\mathfrak{h}}}-\ad_{\xi([x,y]_\g)}\\
                           &=&\ad_{\rho_x(\xi(y))-\rho_y(\xi(x))+[\xi(x),\xi(y)]_{\mathfrak{h}}-\xi([x,y]_\g)}.
\end{eqnarray*}
By $\cen(\h)=0$, we get \eqref{isom3}. Thus, we have \eqref{isom1}-\eqref{isom3}. Therefore, we deduce that $(\mathfrak{g}\oplus \mathfrak{h},[\cdot,\cdot]_{(\rho,\omega)},\phi)$ and $(\mathfrak{g}\oplus \mathfrak{h},[\cdot,\cdot]_{(\rho',\omega')},\phi)$ are isomorphic diagonal non-abelian extensions of $\g$ by $\h$. The proof is finished. \qed

\section{Obstruction of existence of diagonal non-abelian extensions of Hom-Lie algebras}
In this section, we always assume that the following short exact sequences of Hom-Lie algebra
morphisms
\begin{eqnarray*}
0\lon\Inn(\h)\stackrel{}{\lon}\Der(\h)\stackrel{\pi}{\lon} \Out(\h)\lon 0,&&\\
0\lon \cen(\h)\stackrel{}{\lon}\mathfrak{h}\stackrel{\ad}{\lon} \Inn(\mathfrak{h})\lon 0,&&
\end{eqnarray*}
are diagonal non-abelian extensions. Given a Hom-Lie algebra morphism $\bar{\rho}:\g\lon\Out(\h)$, where $\cen(\h)\not=0$, we consider the
obstruction of existence of non-abelian extensions. By choosing a diagonal section $s$ of $\pi:\Der(\h)\lon\Out(\h)$, we can still define
$\rho$ by \eqref{rho} such that \eqref{ext1} hold. Moreover, we can choose a linear map $\omega:\g\wedge\g\lon\h$ such that \eqref{ext2} and \eqref{ext3} hold. Thus,
$(\mathfrak{g}\oplus \mathfrak{h},[\cdot,\cdot]_{(\rho,\omega)},\phi)$ is a diagonal non-abelian extension of $\g$ by $\h$ if and only if
\begin{eqnarray}\label{obstruction}
\rho_{\phi_{\mathfrak{g}}(x)}(\omega(y,z))-\omega([x,y]_{\mathfrak{g}},\phi_{\mathfrak{g}}(z))+c.p.(x,y,z)=0.
\end{eqnarray}
Let $\dM_{\rho}$ be the formal coboundary operator associated to $\rho$. Then we have $$(\dM_{\rho}\omega)(x,y,z)=\rho_{\phi_{\mathfrak{g}}(x)}(\omega(y,z))-\omega([x,y]_{\mathfrak{g}},\phi_{\mathfrak{g}}(z))+c.p.(x,y,z).$$
Therefore, $(\mathfrak{g}\oplus \mathfrak{h},[\cdot,\cdot]_{(\rho,\omega)},\phi)$ is a diagonal non-abelian extension of $\g$ by $\h$ if and only if $\dM_{\rho}\omega=0.$

\begin{defi}
Let $\bar{\rho}:\g\lon\Out(\frkh)$ be a Hom-Lie algebra morphism. We call $\bar{\rho}$ an {\bf extensible homomorphism} if there exists a diagonal section $s$ of $\pi:\Der(\h)\lon\Out(\h)$ and linear map $\omega:\g\wedge\g\lon\h$
such that \eqref{ext1}-\eqref{ext3} and \eqref{obstruction} hold.
\end{defi}

For all $u\in\cen(\h)$, it is
obvious that $\phi_{\mathfrak{h}}(u)\in\cen(\h)$. For $v\in\h$, we have
\begin{eqnarray*}
[\rho_{x}(u),v]_\h&=&[(\Ad_{\phi_{\mathfrak{h}}^{-1}}\rho_{\phi_\g(x)})(u),\phi_{\mathfrak{h}}(\phi_{\mathfrak{h}}^{-1}(v))]_\h\\
                  &=&\rho_{\phi_\g(x)}([u,\phi_{\mathfrak{h}}^{-1}(v)]_\h)-[\phi_{\mathfrak{h}}(u),(\Ad_{\phi_{\mathfrak{h}}^{-1}}\rho_{\phi_\g(x)})(\phi_{\mathfrak{h}}^{-1}(v))]_\h\\
                  &=&0.
\end{eqnarray*}
Thus, we have $\rho_{x}(u)\in\cen(\h)$. Therefore, we can define $\hat{\rho}:\g\lon\gl(\cen(\h))$ by
\begin{eqnarray*}
\hat{\rho}_x\triangleq\rho_x|_{\cen(\h)}.
\end{eqnarray*}
By \eqref{ext1} and \eqref{ext3}, we obtain that $\hat{\rho}$ is a Hom-Lie morphism from $\g$ to $\gl(\cen(\h))$. By Theorem \ref{thm:Hom-Lie rep}, $\hat{\rho}$ is a representation of $(\g,[\cdot,\cdot]_\g,\phi_\g)$ on $\cen(\h)$ with respect to $\phi_{\mathfrak{h}}|_{\cen(\h)}$. By \eqref{ext-isom2}, we deduce that different diagonal sections of $\pi$ give the same representation of $\g$ on $\cen(\h)$ with respect to $\phi_{\mathfrak{h}}|_{\cen(\h)}$. In the
sequel, we always assume that $\hat{\rho}$ is a representation of $\g$ on $\cen(\h)$ with respect to $\phi_{\mathfrak{h}}|_{\cen(\h)}$, which is induced by $\bar{\rho}$. By \eqref{cen}, we have $(\dM_{\rho}\omega)(x,y,x)\in\cen(\h)$. Thus, we have $\dM_{\rho}\omega\in C_{\phi_{\g},\phi_{\mathfrak{h}}\mid_{\cen(\h)}}^3(\g;\cen(\h))$. Moreover, we have the following lemma.

\begin{lem}\label{3-cocycle}
$\dM_{\rho}\omega$ is a $3$-cocycle on $\g$ with the coefficient in $\cen(\h)$ and the cohomology class
$[\dM_{\rho}\omega]$ does not depend on the choices of the diagonal section $s$ of $\pi:\Der(\h)\lon\Out(\h)$ and $\omega$ that we
made.
\end{lem}

\pf For all $x,y,z,t\in\g$, by straightforward computations, we have
\begin{eqnarray*}
&&\dM_{\hat{\rho}}(\dM_{\rho}\omega)(x,y,z,t)\\&=&\rho_{\phi^2_{\g}(x)}(\dM_{\rho}\omega(y,z,t))-\rho_{\phi^2_{\g}(y)}(\dM_{\rho}\omega(x,z,t))+\rho_{\phi^2_{\g}(z)}(\dM_{\rho}\omega(x,y,t))-\rho_{\phi^2_{\g}(t)}(\dM_{\rho}\omega(x,y,z))\\
&&-(\dM_{\rho}\omega)([x,y]_{\g},\phi_{\g}(z),\phi_{\g}(t))+(\dM_{\rho}\omega)([x,z]_{\g},\phi_{\g}(y),\phi_{\g}(t))-(\dM_{\rho}\omega)([x,t]_{\g},\phi_{\g}(y),\phi_{\g}(z))\\
&&-(\dM_{\rho}\omega)([y,z]_{\g},\phi_{\g}(x),\phi_{\g}(t))+(\dM_{\rho}\omega)([y,t]_{\g},\phi_{\g}(x),\phi_{\g}(z))-(\dM_{\rho}\omega)([z,t]_{\g},\phi_{\g}(x),\phi_{\g}(y)).
\end{eqnarray*}
By the definition of $\dM_{\rho}\omega$, we have $60$ terms in the right hand side of above formula. Fortunately, we can cancel the following terms
\begin{eqnarray*}
-\rho_{\phi^2_{\g}(x)}\big(\omega([y,z]_{\g},\phi_{\g}(t))\big)-\rho_{\phi^2_{\g}(x)}\big(\omega(\phi_{\g}(t),[y,z]_{\g})\big)&=&0,\\
\omega([\phi_{\g}(z),\phi_{\g}(t)]_{\g},\phi_{\g}([x,y]_{\g}))+\omega([\phi_{\g}(x),\phi_{\g}(y)]_{\g},\phi_{\g}([z,t]_{\g}))&=&0,\\
\omega([\phi_{\g}(t),[x,y]_{\g}]_{\g},\phi^2_{\g}(z))+\omega([[x,t]_{\g},\phi_{\g}(y)]_{\g},\phi^2_{\g}(z))-\omega([[y,t]_{\g},\phi_{\g}(x)]_{\g},\phi^2_{\g}(z))&=&0.
\end{eqnarray*}
By \eqref{ext1}-\eqref{ext2},  the above formula reduces to the following:
\begin{eqnarray*}
([\rho_{\phi_{\g}(x)},\rho_{\phi_{\g}(y)}]_{\phi_\h}-\rho_{\phi_{\g}([x,y]_{\g})})\big(\omega(\phi_{\g}(z),\phi_{\g}(t))\big)+([\rho_{\phi_{\g}(x)},\rho_{\phi_{\g}(z)}]_{\phi_\h}-\rho_{\phi_{\g}([x,z]_{\g})})\big(\omega(\phi_{\g}(t),\phi_{\g}(y))\big)\\
+([\rho_{\phi_{\g}(x)},\rho_{\phi_{\g}(t)}]_{\phi_\h}-\rho_{\phi_{\g}([x,t]_{\g})})\big(\omega(\phi_{\g}(y),\phi_{\g}(z))\big)+([\rho_{\phi_{\g}(y)},\rho_{\phi_{\g}(z)}]_{\phi_\h}-\rho_{\phi_{\g}([y,z]_{\g})})\big(\omega(\phi_{\g}(x),\phi_{\g}(t))\big)\\
-([\rho_{\phi_{\g}(y)},\rho_{\phi_{\g}(t)}]_{\phi_\h}+\rho_{\phi_{\g}([y,t]_{\g})})\big(\omega(\phi_{\g}(x),\phi_{\g}(z))\big)+([\rho_{\phi_{\g}(z)},\rho_{\phi_{\g}(t)}]_{\phi_\h}-\rho_{\phi_{\g}([z,t]_{\g})})\big(\omega(\phi_{\g}(x),\phi_{\g}(y))\big).
\end{eqnarray*}
Since $[\rho_{x},\rho_{y}]_{\phi_{\mathfrak{h}}}-\rho_{[x,y]_{\mathfrak{g}}}=\ad_{\omega(x,y)}$ and $\phi_{\g}$ is an algebra morphism, this is rewritten as follows:
\begin{eqnarray*}
&&[\omega(\phi_{\g}(x),\phi_{\g}(y)),\omega(\phi_{\g}(z),\phi_{\g}(t))]_{\mathfrak{h}}+[\omega(\phi_{\g}(x),\phi_{\g}(z)),\omega(\phi_{\g}(t),\phi_{\g}(y))]_{\mathfrak{h}}\\
&&+[\omega(\phi_{\g}(x),\phi_{\g}(t)),\omega(\phi_{\g}(y),\phi_{\g}(z))]_{\mathfrak{h}}+[\omega(\phi_{\g}(y),\phi_{\g}(z)),\omega(\phi_{\g}(x),\phi_{\g}(t))]_{\mathfrak{h}}\\
&&-[\omega(\phi_{\g}(y),\phi_{\g}(t)),\omega(\phi_{\g}(x),\phi_{\g}(z))]_{\mathfrak{h}}+[\omega(\phi_{\g}(z),\phi_{\g}(t)),\omega(\phi_{\g}(x),\phi_{\g}(y))]_{\mathfrak{h}}\\
&&=0.
\end{eqnarray*}
Thus, we obtain $\dM_{\rho}\omega\in \huaZ^3(\g;\hat{\rho})$.

Now Let us check that the cohomology class $[\dM_{\rho}\omega]$ does not depend on the choices of the diagonal section $s$ of $\pi:\Der(\h)\lon\Out(\h)$ and $\omega$ that we made. Let $s'$ be another diagonal section of $\pi$, we have $\rho'_x\in\Der(\h)$ and choose $\omega'$, such that \eqref{ext1}-\eqref{ext3} hold. We are going to prove that $[\dM_{\rho}\omega]=[\dM_{\rho'}\omega']$. Since $s$ and $s'$ are diagonal sections of $\pi$, we have linear map $b:\g\lon\h$ such that
$$
b\circ\phi_{\g}=\phi_{\h}\circ b ,\,\,\,\,\rho'_x=\rho_x+\ad_{b(x)}.
$$
We define $\omega^*$ by
$$
\omega^*(x,y)=\omega(x,y)+\rho_x(b(y))-\rho_y(b(x))+[b(x),b(y)]_{\h}-b[x,y]_{\g}.
$$
By straightforward computations, we obtain that \eqref{ext2} and \eqref{ext3} hold for $\rho',\omega^*$. For all $x,y,z\in\g$, we have
\begin{eqnarray*}
&&(\dM_{\rho'}\omega^*-\dM_{\rho}\omega)(x,y,z)\\&=&(\rho'_{\phi_{\g}(x)}\circ\omega^*(y,z)-\omega^*([x,y]_{\g},\phi_{\g}(z))-(\rho_{\phi_{\g}(x)}\circ\omega(y,z)-\omega([x,y]_{\g},\phi_{\g}(z)))+c.p.\\
                                             &=&A(x,y,z)+A(y,z,x)+A(z,x,y)+B(x,y,z)+B(y,z,x)+B(z,x,y)+C(x,y,z),
\end{eqnarray*}
where
\begin{eqnarray*}
A(x,y,z)&=&\rho_{\phi_{\g}(x)}\circ\rho_y(b(z))-\rho_{\phi_{\g}(y)}\circ\rho_x(b(z))-\rho_{[x,y]_{\g}}(b(\phi_{\g}(z)))+[b(\phi_{\g}(z)),\omega(x,y)]_{\h},\\
B(x,y,z)&=&\rho_{\phi_{\g}(x)}([b(y),b(z)]_{\h})-[\rho_x(b(y)),b(\phi_{\g}(z))]_{\h}-[b(\phi_{\g}(y)),\rho_x(b(z))]_{\h},\\
C(x,y,z)&=&[b(\phi_{\g}(x)),[b(y),b(z)]_{\h}]_{\h}+[b(\phi_{\g}(y)),[b(z),b(x)]_{\h}]_{\h}+[b(\phi_{\g}(z)),[b(x),b(y)]_{\h}]_{\h}\\
&&+b([[x,y]_{\g},\phi_{\g}(z)]_{\g})+b([[y,z]_{\g},\phi_{\g}(x)]_{\g})+b([[z,x]_{\g},\phi_{\g}(y)]_{\g}).
\end{eqnarray*}
By \eqref{ext1} and \eqref{ext3}, we have $A(x,y,z)=0$. Since $\rho_{\phi_{\g}(x)}$ is a derivation, we obtain $B(x,y,z)=0$. Since $b\circ\phi_{\g}=\phi_{\h}\circ b$ and $\g,\h$ are Hom-Lie algebras, we get $C(x,y,z)=0$. Thus, we have $\dM_{\rho}\omega=\dM_{\rho'}\omega^*$. Since the equations
\eqref{ext2} and \eqref{ext3} hold for $\rho',\omega^*$ and $\rho',\omega'$ respectively, we have
$$
[\rho'_x,\rho'_y]_{\phi_{\h}}-\rho'_{[x,y]_{\g}}=\ad_{\omega'(x,y)}=\ad_{\omega^*(x,y)}.
$$
Thus, we have $\ad_{(\omega'-\omega^*)(x,y)}=0$. Moreover, we have $(\omega'-\omega^*)(x,y)\in\cen(\h)$. By \eqref{ext2}, we can define
$\tau\in C_{\phi_{\g},\phi_{\mathfrak{h}}\mid_{\cen(\h)}}^2(\g;\cen(\h))$ by
$$
\tau(x,y)=\omega'(x,y)-\omega^*(x,y).
$$
Thus, we have
$$
\dM_{\rho'}\omega'-\dM_{\rho}\omega=\dM_{\rho'}\omega'-\dM_{\rho'}\omega^*=\dM_{\rho'}(\omega'-\omega^*)=\dM_{\rho'}\tau=\dM_{\hat{\rho}}\tau.
$$
Therefore, we obtain $[\dM_{\rho}\omega]=[\dM_{\rho'}\omega']$. The proof is finished. \qed\vspace{3mm}

Now we are ready to give the main result in this paper, namely the obstruction of a Hom-Lie algebra homomorphism $\bar{\rho}:\g\lon\Out(\frkh)$ being extensible is given by the cohomology class $[\dM_{\rho}\omega]\in \huaH^3(\g;\hat{\rho})$.

\begin{thm}
Let $\bar{\rho}:\g\lon\Out(\frkh)$ be a Hom-Lie algebra morphism. Then $\bar{\rho}$ is an
extensible homomorphism if and only if
$$[\dM_{\rho}\omega]=[0].$$
\end{thm}

\pf
Let $\bar{\rho}:\g\lon\Out(\frkh)$ be an extensible Hom-Lie algebra morphism. Then we can
choose a diagonal section $s$ of $\pi:\Der(\h)\lon\Out(\h)$ and define $\rho$ by \eqref{rho}. Moreover, we can choose a linear map $\omega:\g\wedge\g\lon\h$ such that \eqref{ext1}-\eqref{ext3} hold. Since $\bar{\rho}$ is extensible, we have $\dM_{\rho}\omega=0$,  which
implies that $[\dM_{\rho}\omega]=[0]$.

Conversely, if $[\dM_{\rho}\omega]=[0]$, then there exists $\sigma\in C_{\phi_{\g},\phi_{\mathfrak{h}}\mid_{\cen(\h)}}^2(\g;\cen(\h))$ such that $\dM_{\rho}\omega=\dM_{\hat{\rho}}\sigma.$
Thus, we have
\begin{eqnarray*}
\dM_{\rho}(\omega-\sigma)=\dM_{\rho}\omega-\dM_{\rho}\sigma=\dM_{\rho}\omega-\dM_{\hat{\rho}}\sigma=0.
\end{eqnarray*}
Since $\sigma\in C_{\phi_{\g},\phi_{\mathfrak{h}}\mid_{\cen(\h)}}^2(\g;\cen(\h))$, we also have
\begin{eqnarray*}
\label{}\phi_{\mathfrak{h}}\circ (\omega-\sigma)&=&(\omega-\sigma)\circ\phi_{\mathfrak{g}}^{\otimes 2},\\
\label{}[\rho_{x},\rho_{y}]_{\phi_{\mathfrak{h}}}-\rho_{[x,y]_{\g}}&=&\ad_{(\omega-\sigma)(x,y)}.
\end{eqnarray*}
By Proposition \ref{extension}, we can construct a Hom-Lie algebra $(\mathfrak{g}\oplus \mathfrak{h},[\cdot,\cdot]_{(\rho,\omega-\sigma)},\phi)$. Therefore, $\bar{\rho}$ is an extensible morphism. The proof
is finished. \qed\vspace{3mm}

The following theorem classifies diagonal non-abelian extensions of $\g$ by $\h$ once they exist.

\begin{thm}
Let $\bar{\rho}:\g\lon\Out(\frkh)$ be an extensible morphism. Then isomorphism classes
of diagonal non-abelian extensions of $\g$ by $\h$ induced by $\bar{\rho}$ are parameterized by $\huaH^2(\g;\hat{\rho})$.
\end{thm}

\pf
Since $\bar{\rho}$ is an extensible morphism, we can choose a diagonal section $s$ of $\pi$ and define $\rho$ by
\eqref{rho}. We choose a linear map $\omega$
such that \eqref{ext1}-\eqref{ext3} hold and $\dM_{\rho}\omega=0$. Thus, the Hom-Lie algebra $(\mathfrak{g}\oplus \mathfrak{h},[\cdot,\cdot]_{(\rho,\omega)},\phi)$ define by \eqref{dbr} and \eqref{dmo} is a diagonal non-abelian extension of $\g$ by $\frkh$, which is induced by $\bar{\rho}$.  Let $s'$ be another diagonal section of $\pi$ and define $\rho'$ by \eqref{rho}. We also choose a linear map $\omega'$ such that \eqref{ext1}-\eqref{ext3} hold and $\dM_{\rho'}\omega'=0$. Since $s$ and $s'$ are diagonal sections of $\pi$, there exists a linear  map $b:\g\lon\h$ such that
$$
\phi_{\h}\circ b=b\circ\phi_{\g},\,\,\,\,
\rho_x =\rho'_x+\ad_{b(x)}.
$$
We define $\omega^*$ by
$$
\omega^*(x,y)=\omega'(x,y)+\rho'_x(b(y))-\rho'_y(b(x))+[b(x),b(y)]_{\h}-b[x,y]_{\g}.
$$
By the computation in Lemma \ref{3-cocycle}, we have $\dM_{\rho}\omega^*=\dM_{\rho'}\omega'=0$. Thus, the
Hom-Lie algebra $(\mathfrak{g}\oplus \mathfrak{h},[\cdot,\cdot]_{(\rho,\omega^*)},\phi)$ constructed from $\rho,\omega^*$ is isomorphic to the Hom-Lie algebra $(\mathfrak{g}\oplus \mathfrak{h},[\cdot,\cdot]_{(\rho',\omega,)},\phi)$ constructed from $\rho',\omega'$. Thus, we only need to study the Hom-Lie algebras constructed from a fix diagonal section $s$. For all $\tilde{\omega}$, which satisfy \eqref{ext2}-\eqref{ext3} and $\dM_{\rho}\tilde{\omega}=0$, we  define
$$\lambda=\omega-\tilde{\omega}\in C_{\phi_{\g},\phi_{\mathfrak{h}}\mid_{\cen(\h)}}^2(\g;\cen(\h)).$$
Moreover, we have
$$\dM_{\hat{\rho}}\lambda=\dM_{\rho}(\omega-\tilde{\omega})=0-0=0.$$
which implies that $\omega-\tilde{\omega}\in\huaZ^2(\g;\hat{\rho})$.

Moreover, if the Hom-Lie algebra $(\mathfrak{g}\oplus \mathfrak{h},[\cdot,\cdot]_{(\rho,\omega)},\phi)$  constructed from $\rho,\omega$ is isomorphic to the Hom-Lie algebra $(\mathfrak{g}\oplus \mathfrak{h},[\cdot,\cdot]_{(\rho,\tilde{\omega})},\phi)$ constructed from $\rho,\tilde{\omega}$. Then there exists a linear map $b:\g\lon\frkh$ which does not change $\rho$, i.e. $b:\g\lon\cen(\h)$, such that
\begin{eqnarray*}
 \label{}\phi_{\h}\circ b&=&b\circ\phi_{\g},\\
  %\label{exiso2}\rho_x -\rho_x&=&\ad_{b(x)}, \\
\label{} \omega-\tilde{\omega}&=&\dM_{\rho}b.
\end{eqnarray*}
This is equivalent to that $\omega-\tilde{\omega}\in\huaB^2(\g;\hat{\rho}).$ Thus, isomorphism classes of diagonal non-abelian extensions of $\g$ by $\h$ induced by $\bar{\rho}$ are parameterized by
$\huaH^2(\g;\hat{\rho})$.  \qed

\begin{cor}
The isomorphism classes of diagonal non-abelian extensions of a Hom-Lie algebra $\g$ by
a Hom-Lie algebra $\h$ correspond bijectively to the set of pairs $(\bar{\rho},[\kappa])$, where $\bar{\rho}$ is an extensible
morphism from $\g$ to $\Out(\h)$ and $[\kappa]\in\huaH^2(\g;\hat{\rho})$.
\end{cor}

\end{document}